\documentclass[11pt]{amsart}
\usepackage[utf8]{inputenc}
\pdfoutput=1
\usepackage{navigator}

\usepackage{enumitem}
\usepackage[T1]{fontenc}

\usepackage{hyperref}
\usepackage[english]{babel}

\usepackage{amssymb,amsfonts,stmaryrd,mathrsfs}
\usepackage{amsmath,latexsym,amsthm,mathtools,bbm}
\usepackage{tikz-cd}
\usepackage{relsize,exscale}
\usepackage{color}
 
\usepackage{graphicx}
\usepackage[export]{adjustbox}
\usepackage{pgf,tikz} 
\usepackage{mathrsfs}
\usepackage{bm}
\usepackage{colortbl}
\usetikzlibrary{arrows}
\usepackage[capitalize]{cleveref}
\usepackage{graphicx}
\usepackage{multicol}

\newtheorem{thm}{Theorem}[section]
\newtheorem{cor}[thm]{Corollary}
\newtheorem{lem}[thm]{Lemma}
\newtheorem{prop}[thm]{Proposition}
\newtheorem*{conj_MJ}{The Matching-Jack conjecture}
\newtheorem*{conj_MJ2}{The Matching-Jack conjecture (combinatorial formulation)}
\newtheorem*{b-conj}{The $b$-conjecture}

\theoremstyle{definition}
\newtheorem{defi}{Definition}[section]

\theoremstyle{remark}
\newtheorem{rmq}{Remark}
\newtheorem{exe}{Example}[section]

\definecolor{uibred}  {HTML}{db3f3d}
\definecolor{uibblue} {HTML}{4ea0b7}
\definecolor{uibgreen}{HTML}{274e13}
\definecolor{uibgray} {HTML}{d0cac2}
\definecolor{uiblink} {HTML}{00769E}

\definecolor{orange} {HTML}{ff8c1a}
\definecolor{vert} {HTML}{009933}
\definecolor{bleu} {HTML}{0033cc}

\author{Houcine Ben Dali}
\date{}

\DeclareMathOperator{\rk}{rk}
\DeclareMathOperator{\Z}{\mathcal Z}
\DeclareMathOperator{\Span}{Span}
\DeclareMathOperator{\tc}{t}
\DeclareMathOperator{\cf}{\mathfrak c}
\DeclareMathOperator{\mf}{\mathfrak m}
\DeclareMathOperator{\ff}{\mathfrak f}
\DeclareMathOperator{\gf}{\mathfrak g}
\DeclareMathOperator{\F}{\mathfrak F}
\DeclareMathOperator{\mc}{ \mathit{\overline{ c}}}

 \address{Université de Lorraine\\ CNRS\\ IECL\\ F-54000 Nancy\\ France} 
\address{Université de Paris\\ CNRS\\ IRIF\\ F-75006\\ Paris\\ France}

\email{houcine.ben-dali@univ-lorraine.fr}

\begin{document}

\title[Integrality in the Matching-Jack conjecture]{Integrality in the Matching-Jack conjecture  and the Farahat-Higman algebra}
\begin{abstract}
Using Jack polynomials, Goulden and Jackson have introduced a  one parameter deformation $\tau_b$ of the generating series of bipartite maps, which generalizes the partition function of $\beta$-ensembles of random matrices. The Matching-Jack conjecture suggests that the coefficients $c^\lambda_{\mu,\nu}$ of the function $\tau_b$ in the power-sum basis are non-negative integer polynomials in the deformation parameter $b$. 
Do\l{}{\k{e}}ga and Féray have proved in 2016 the "polynomiality" part in the Matching-Jack conjecture, namely that coefficients $c^\lambda_{\mu,\nu}$ are in $\mathbb{Q}[b]$. 
In this paper, we prove the "integrality" part, i.e that the coefficients $c^\lambda_{\mu,\nu}$ are in $\mathbb{Z}[b]$.

The proof is based on a recent work of the author that deduces the Matching-Jack conjecture for marginal sums $\mc^\lambda_{\mu,l}$ from an analog result for the $b$-conjecture, established in 2020 by Chapuy and Do\l{}{\k{e}}ga.  
A key step in the proof involves a new connection with the graded Farahat-Higman algebra.    
\end{abstract}

\maketitle

\section{Introduction:}
\renewcommand{\labelitemi}{\textbullet}

\subsection{Coefficients \texorpdfstring{$c^\lambda_{\mu,\nu}$}{} and the Matching-Jack conjecture.}\label{ssec coef c}
Jack polynomials $J^{(\alpha)}_\theta$ are symmetric functions that depend on a deformation parameter $\alpha$ (see \cite{J,Stan89,Mac}). 
We consider the power-sum functions  $\mathbf{p}:=(p_1,p_2,..)$, $\mathbf{q}:=(q_1,q_2,..)$ and $\mathbf{r}:=(r_1,r_2,..)$ associated to three different alphabets, and we denote by $J^{(\alpha)}_\theta$ the Jack polynomial of parameter $\alpha$ expressed in the power-sum basis.
In \cite{GJ96b}, Goulden and Jackson have introduced a function $\tau_b$ which depends on the parameter $b:=\alpha-1$ and defined by:
$$\tau_b(t,\mathbf{p},\mathbf{q},\mathbf{r}):=\sum_{n\geq0}t^n\sum_{\theta\vdash n}\frac{1}{j^{(\alpha)}_\theta}J^{(\alpha)}_\theta(\mathbf{p})J^{(\alpha)}_\theta(\mathbf{q})J^{(\alpha)}_\theta(\mathbf{r}),$$ 
where the second sum is taken over integer partitions of size $n$, and
where  $j^{(\alpha)}_\theta$ is the square norm of the Jack polynomial $J_\theta^{(\alpha)}$, see \cref{sec partitions}.
The function $\tau_b$ has been introduced as an interpolation between the generating series of bipartite maps on orientable surfaces and general surfaces, obtained by setting respectively the parameter $b$ to 0 and 1. 
This function is also related to the $\beta$-ensembles of random matrix theory, which can be obtained by some specializations of the variables $\mathbf{q}$ and $\mathbf{r}$, see \cite{La09,CD20}. Under some specializations, the function $\tau_b$ satisfies a family of partial differential equations known in theoretical physics as the Virasoro constraints, see \cite{KZ15,BCD21}. For $b=0$, the function $\tau_b$ is related to the generating series of weighted Hurwitz numbers \cite{HGP17}, with strong links to the topological recursion \cite{ACEH20,B20}.

For a partition $\lambda$, we write $p_\lambda:=p_{\lambda_1}p_{\lambda_2}\dots$.
We use the same notation for $\mathbf{q}$ and $\mathbf{r}$.
The coefficients $c^\lambda_{\mu,\nu}(b)$ are defined as the coefficients of the function $\tau_b$ in the power-sum basis:
\begin{equation}\label{eqtau}
\tau_b(t,\mathbf{p},\mathbf{q},\mathbf{r})=1+\sum_{n\geq 1}t^n\sum_{\lambda,\mu,\nu\vdash n}\frac{c^\lambda_{\mu,\nu}(b)}{z_\lambda(1+b)^{\ell(\lambda)}}p_\lambda q_\mu r_\nu,
\end{equation}
where $\ell(\lambda)$ denotes the length of $\lambda$ and $z_\lambda$ is a normalizing factor defined in \cref{sec partitions}.
The coefficients $c^\lambda_{\mu,\nu}(b)$ are the main objects of the Matching-Jack conjecture, formulated by Goulden and Jackson \cite[Conjecture 3.5]{GJ96b}.
\begin{conj_MJ}[\cite{GJ96b}]
For all partitions $\lambda,\mu,\nu$ of size $n\geq1$, the coefficient $c^\lambda_{\mu,\nu}$ is a polynomial in $b$ with non-negative integer coefficients.
\end{conj_MJ}
The Matching-Jack conjecture is equivalent to saying that $\tau_b$ is a generating series of matchings with some particular weights (see \cref{ssec matchings}). 
There exists a "connected" version of this conjecture called the \textit{$b$-conjecture.} The $b$-conjecture suggests that the function $(1+b)\frac{t\partial}{\partial t}\log(\tau_b)$ has also a positivity property and that it enumerates bipartite maps with some weights \cite[Conjectures 6.2 and 6.3]{GJ96b}.  Since bipartite maps can be encoded with matchings (see \textit{e.g} \cite{GJ96a,DFS14,BD21}), the Matching-Jack conjecture and the $b$-conjecture are related. However, no implication between them has been proved and both of them are still open. See \cref{ssec bconj} for further details about the connection between the two conjectures.

\subsection{Former results and main theorem.}
In addition to the special cases $b=0$ and $b=1$ that follow from connections with representation theory respectively of $\mathfrak{S}_n$ and the Gelfand pair $(\mathfrak{S}_{2n},H_n)$ \cite{LZ04,GJ96b}, several partial results related to the Matching-Jack conjecture have been established; we cite here some of them. 

It follows from the definitions that the coefficients $c^\lambda_{\mu,\nu}$ are rational functions in $b$. The polynomiality is however not immediate and has been proved by Féray and Do\l{}{\k{e}}ga.

\begin{thm}\cite[Corollary 4.2]{DF16}\label{thm polynomiality}
For every $n\geq1$, and every $\lambda,\mu,\nu\vdash n$, the coefficient $c^\lambda_{\mu,\nu}$ is a polynomial in $b$ with rational coefficients. 
\end{thm}

Our main result establishes the integrality  of the coefficients of these polynomials.

\begin{thm}[Main result]\label{thm integrality}
For every $n\geq 1$, and every $\lambda,\mu,\nu\vdash n$, the coefficient $c^\lambda_{\mu,\nu}$ is a polynomial in $b$ with integer coefficients. 
\end{thm}

 Since the approach used here is independent from the one considered in \cite{DF16}, it gives a new proof of \Cref{thm polynomiality}. Unfortunately, the non-negativity of the coefficients of $c^\lambda_{\mu,\nu}$ as polynomials in $b$ --the remaining part of the Matching-Jack conjecture-- seems to be out of reach with our approach and requires new ideas.

Although we are first to establish the integrality in the Matching-Jack conjecture in full generality, integrality (and in fact, non-negativity) have already been proved for certain marginal sums of the coefficients $c^\lambda_{\mu,\nu}$ (corresponding to some specializations of the function $\tau_b$). This result has been established by Chapuy and Do\l{}{\k{e}}ga in \cite{CD20} in the case of the $b$-conjecture, and an analog for the Matching-Jack conjecture has been deduced in \cite{BD21}.

  \begin{thm}\cite[Theorem 1.9]{BD21}\label{thm mar}
    For every $n,l\geq1$ and for every $\lambda,\mu\vdash n$, the marginal coefficient
    $\mc^\lambda_{\mu,l}$ defined by 
    $$\mc^\lambda_{\mu,l}:=\sum_{\ell(\nu)=l}c^\lambda_{\mu,\nu}\text{ ,}$$
    is a polynomial in $b$ with non-negative integer coefficients.
    \end{thm}
The previous theorem covers other partial results in this direction (see \cite{KV16,KPV18}).
Our proof of \cref{thm integrality} strongly relies on \cref{thm mar}. In fact, we deduce the integrality of the coefficients $c^\lambda_{\mu,\nu}$ from the integrality of their marginal sums $\mc^\lambda_{\mu,l}$, thanks to a new connection between them and the Farahat-Higman algebra. This two-step method is a key feature of our proof; the approach used to establish integrality for the marginal coefficients $\mc^\lambda_{\mu,l}$ relies on the existence of functional equations for specializations of the function $\tau_b$ (Tutte-like equations), which are unknown in the case of the function $\tau_b$ itself.

\subsection{The Farahat-Higman algebra}
The Farahat-Higman algebra was introduced in \cite{FH59} in order to study the structure coefficients of the conjugacy classes $C_\mu(n)$ in the center of the symmetric group algebra $Z(\mathbb{Z}\mathfrak{S}_n)$.
It has been shown that the Farahat-Higman algebra is isomorphic to the algebra of integral symmetric functions, see \cite{GJ94,CGS04}. It is also related to the algebra of partial permutations introduced by Ivanov and Kerov in \cite{IK99}.

In \Cref{sec FHA}, we will consider a graded version of the Farahat-Higman algebra that we denote $\mathcal{Z}_\infty$ and that has been introduced in   \cite[Example 24, page 131]{Mac}.
In \Cref{thm g}, we exhibit a new basis of $\Z_\infty$, which is useful in the proof of the main theorem.
Since the Farahat-Higman is of independent interest, \cref{thm g} might reveal itself useful in the future, independently of its application in the current paper.

\subsection{Steps of the proof.}\label{ssec steps of the proof}
A key tool of the proof of \Cref{thm integrality} is the following multiplicativity property which is a consequence of the orthogonality of Jack polynomials (see also \Cref{prop multiplicativity}); for all partitions $\lambda,\mu,\nu$ of the same size $n\geq 1$ and for every $l\geq1$, one has
\begin{equation}\label{eq multiplicativity}
  \sum_{\kappa\vdash n}c^\lambda_{\mu,\kappa}\mc^\kappa_{\nu,l}=\sum_{\theta\vdash n} \mc^\lambda_{\theta,l} c^\theta_{\mu,\nu}.
\end{equation}
 The equations \eqref{eq multiplicativity} will be considered as a system of linear equations that allows to recover $c^\lambda_{\mu,\nu}$ from the marginal coefficients $\mc^\lambda_{\mu,l}$ (see \Cref{rmq eqsystem}).
 We now give the key steps of the proof of \Cref{thm integrality}.

\begin{itemize}

     \item We prove that for a particular choice of instances of \Cref{eq multiplicativity} (i.e for a subset of quadruples $(\lambda,\mu,\nu,l)$), we obtain a square linear system. The matrix of this system is block triangular. 
    \item We prove that the diagonal blocks of the matrix encoding this system, denoted $\mathcal{Q}^{(r)}$, contain some coefficients $\tc^\rho_\pi$, that are independent from $b$ (see \Cref{sec top coef}).
    \item We prove that the matrix $\mathcal{Q}^{(r)}$ is invertible in $\mathbb Z$ by seeing it as a change-of-basis matrix in the graded Farahat-Higman algebra $\Z_\infty$ (see \Cref{prop Q tr mat} and \Cref{thm g}). 
\end{itemize}

\subsection{Multiplicativity property for matchings and some consequences of the main result}
As a consequence of our techniques, we obtain a new proof for the "Matching-Jack conjecture"\footnote{We recall that the cases $b=0$ and $b=1$ have preceded the conjecture \cite{GJ96a}, hence the quotes.} for $b=0$ and $b=1$. Unlike the proof given in \cite{GJ96a}, the proof we give here does not use representation theory\footnote{A more intricate representation-free proof can also be obtained using the arguments provided in \cite{CD20} (private communication with Guillaume Chapuy and Maciej Do\l{}{\k{e}}ga).}. This new proof, detailed in \Cref{sec b=0 and b=1}, relies on the following three ingredients: \begin{itemize}
    \item We observe that, at a combinatorial level, matchings satisfy a multiplicativity property of the same form of \Cref{eq multiplicativity}.
    \item  We use the combinatorial interpretation of the coefficients $\mc^\lambda_{\mu,l}$ in terms of matchings given in \cite{BD21}, see \cref{thm marg b=0 and b=1}.
    \item As in the proof of the main result, we use the fact that \Cref{eq multiplicativity} entirely determines the coefficients $c^\lambda_{\mu,\nu}$ from their marginal sums, see \cref{prop multiplicativity}.
\end{itemize}

As explained above, the Matching-Jack conjecture is related to the $b$-conjecture. Unfortunately,
we were not able to use the approach detailed in \cref{ssec steps of the proof} to obtain integrality in the $b$-conjecture. However, it is possible to obtain the integrality for the "cumulants" of $c^\lambda_{\mu,\nu}$, which up to rescaling by a factor of the form $\frac{n}{z_\lambda(1+b)^{\ell(\lambda)-1}}$, give the coefficients of the $b$-conjecture (see \cref{ssec bconj}).

\subsection{Outline of the paper.}
In \Cref{sec top coef}, we introduce the top coefficients $\tc^\rho_\pi$ and we formulate \Cref{thm Qinv} which is a key step in the proof of the main result.
In \cref{sec proof main result} we prove that \Cref{thm Qinv} implies \Cref{thm integrality}.  We consider a graded version of the Farahat-Higman algebra in \Cref{sec FHA} and we use it to prove \Cref{thm Qinv}. In \cref{sec b=0 and b=1}, we give a combinatorial interpretation for the multiplicativity property when $b=0$ and $b=1$ and we use it to obtain a new proof for the Matching-Jack conjecture in these cases. In \cref{sec cumulants}, we give some consequences of the main result.

\section{Preliminaries}
\subsection{Some notation}\label{sec partitions}
For this subsection, we refer the reader to \cite{Stan89,Mac} for further details. An \textit{integer partition} $\lambda=[\lambda_1,\lambda_2,...,\lambda_l]$ is a sequence of weakly decreasing positive integers $\lambda_1\geq\lambda_2\geq ...\geq\lambda_l>0$. The integer $l$ is the \textit{length} of the partition $\lambda$, \textit{i.e} its number of parts. It is denoted by $\ell(\lambda)$. The size of $\lambda$ is the integer $|\lambda|:=\lambda_1+\lambda_2+...+\lambda_\ell.$ If $n$ is the \textit{size} of $\lambda$, we say that $\lambda$ is a partition of $n$ and we write $\lambda\vdash n$.  We define the \textit{rank} of a partition $\lambda$ by $\rk(\lambda):=|\lambda|-\ell(\lambda)$. 
We denote by  $\lambda-\mathbf{1}$ the partition  $[\lambda_1-1,\lambda_2-1,\dots]$.
We set for every partition $\lambda$
$$z_\lambda:=\prod_{i\geq1}m_i(\lambda)!i^{m_i(\lambda)},$$
where $m_i(\lambda)$ denotes the number of parts equal to $i$ in $\lambda$.

Let $\lambda$ and $\mu$ be two partitions. We define their \textit{entry-wise sum} by
$$\lambda\oplus \mu=[\lambda_1+\mu_1,\lambda_2+\mu_2,\dots].$$
We also define the \textit{union} partition $\lambda\cup\mu$ as the partition whose parts are obtained by taking the union of the parts of $\lambda$ and $\mu$. In other words, for every $i\geq1$, one has 
$$m_i(\lambda\cup\mu)=m_i(\lambda)+m_i(\mu).$$


\noindent We identify a partition  $\lambda$ with its \textit{Young diagram}, defined by 
$$\lambda:=\{(i,j),1\leq i\leq \ell(\lambda),1\leq j\leq \lambda_i\}.$$
\textit{The conjugate partition} of $\lambda$, denoted $\lambda'$, is the partition associated to the Young diagram obtained by reflecting the diagram of $\lambda$ with respect to the line $j=i$:
$$\lambda':=\{(i,j),1\leq j\leq \ell(\lambda),1\leq i\leq \lambda_j\}.$$

Let $\prec$ denote the \textit{dominance} order, the partial order on partitions of the same size, defined by:
$$\mu\prec \lambda \Longleftrightarrow \mu_1+\mu_2+...+\mu_i\leq \lambda_1+\lambda_2+...+\lambda_i, \forall i\geq1,$$
where we use the convention; $\xi_i=0$ if $i>\ell(\xi)$, for any partition $\xi$.
\begin{defi}\label{def orders}
For every $n\geq0$, we consider a total order $\leq$ on the partitions of size $n$, with the two following properties:
\begin{enumerate}
    \item If $\mu\prec\lambda$, then $\mu\leq \lambda$.
    \item If $\ell(\lambda)<\ell(\mu)$ then $\mu\leq \lambda$.
\end{enumerate}
Such order is well defined since  $\mu\prec\lambda$ implies $\ell(\lambda)\leq \ell(\mu).$

For any partition $\lambda$, we denote by $\lambda'$ the conjugate partition of $\lambda$. We define the dual order $\leq'$ as the total order given by 
$$\mu\leq' \lambda \iff \lambda'\leq\mu'.$$
In particular $\lambda_1 < \mu_1$ implies $\lambda\leq'\mu$.
\end{defi}

Everywhere in this paper, $\alpha$ and $b$ are two parameters related by $\alpha=b+1$. Let $\langle.,.\rangle_\alpha$ denote the scalar product on the ring of symmetric functions, defined by 
$$\langle p_\lambda, p_\mu\rangle_\alpha=\delta_{\lambda,\mu}z_\lambda \alpha^{\ell(\lambda)},$$
where $p_\lambda$ denotes the power-sum function associated to $\lambda$ and $\delta_{\lambda,\mu}$ is the Kronecker delta.

For any partition $\theta$, we denote by $J_\theta^{(\alpha)}$ the Jack polynomial of parameter $\alpha$ associated to the partition $\theta$. Jack polynomials form an orthogonal basis of the ring of symmetric functions, we denote by $j^{(\alpha)}_\theta$ their squared norm:
$$j^{(\alpha)}_\theta:=\langle J^{(\alpha)}_\theta, J^{(\alpha)}_\theta\rangle_\alpha.$$
The quantities $j^{(\alpha)}_\theta$ have an explicit combinatorial expression (see \cite[Theorem 5.8]{Stan89}).

\subsection{Top coefficients \texorpdfstring{$\tc^\rho_\pi$}{}}\label{sec top coef}

As announced in the introduction, the proof of \Cref{thm integrality} involves the resolution of a linear system satisfied by the coefficients $c^\lambda_{\mu,\nu}$. In this section, we introduce the matrices $\mathcal{Q}^{(r)}$ that encode this system. 

We recall that \cref{thm mar} implies that the coefficients $\mc^\lambda_{\mu,l}$ are polynomials in the parameter $b$. The following lemma gives an upper bound on the degree of these polynomials.
\begin{lem}\label{lem bounddeg}
For any partitions $\lambda,\mu\vdash n\geq 1$ and $l\geq 1$, we have the following bound on the degree of $\mc^\lambda_{\mu,l}$ as a polynomial in $b$:
$$\deg_b(\mc^\lambda_{\mu,l})\leq n-l+\ell(\lambda)-\ell(\mu).$$
\begin{proof}
For any partitions $\lambda$ and $\nu$ of size $n$, we have that (see \cite[Lemma 3.2]{GJ96b})
$$\sum_{\mu\vdash n}c^\lambda_{\mu,\nu}=\frac{n!}{z_\nu}(1+b)^{n-\ell(\nu)}.$$
We take the sum over the partitions $\nu$ of length $l$:
$$\sum_{\mu\vdash n}\mc^\lambda_{\mu,l}=\sum_{\ell(\nu)=l}\frac{n!}{z_\nu}(1+b)^{n-\ell(\nu)}.$$
Since the polynomials $\mc^\lambda_{\mu,l}$ have non-negative coefficients (see \Cref{thm mar}), we deduce that $\deg_b\left(\mc^\lambda_{\mu,l}\right)\leq n-l$. On the other hand, from the definition of the coefficients $c^\lambda_{\mu,\nu}$, we see that
$$\frac{c^\lambda_{\mu,\nu}}{z_\lambda(1+b)^{\ell(\lambda)}}=\frac{c^\mu_{\lambda,\nu}}{z_\mu(1+b)^{\ell(\mu)}},$$ hence 
$\deg_b\left(\mc^\lambda_{\mu,l}\right)=\deg_b\left(\frac{z_\lambda(1+b)^{\ell(\lambda)}}{z_\mu(1+b)^{\ell(\mu)}}\mc^\mu_{\lambda,l}\right)\leq n-l+\ell(\lambda)-\ell(\mu).$
\end{proof}
\end{lem}

In this section, we are interested in the coefficients $\mc^\lambda_{\mu,l}$ for which the bound given in \cref{lem bounddeg} is zero. The following lemma gives a stability property of these coefficients when adding parts of size 1. The proof is postponed to \cref{sec FHA}.

\begin{lem}\label{lem add part1}
Let $\kappa,\nu,\theta\vdash m\geq1,$ such that $\rk(\kappa)=\rk(\nu)+\rk(\theta)$. 
Then \begin{enumerate}[label={\normalfont(}\arabic*\normalfont)]
    \item  For every  $n\geq 1$, we have
    $$ c^\kappa_{\nu,\theta}(0)=c^{\kappa\cup1^n}_{\nu\cup{1^n},\theta\cup1^n}(0).$$
    \item If $c^\kappa_{\nu,\theta}(0)\neq 0$, then $m_1(\kappa)\leq m_1(\theta)$, where $m_1(\cdot)$ denotes the number of parts equal~to~1.
\end{enumerate}
\end{lem}

\begin{rmq}
Actually, \cref{lem add part1} holds without the specialization at $b=0$, since under the condition $\rk(\kappa)=\rk(\nu)+\rk(\theta)$, the coefficient $c^\kappa_{\nu,\theta}$ is independent from $b$. This is a consequence of a generalized version of \cref{lem bounddeg} (see \cite[Corollary 4.2]{DF16}).
\end{rmq}
We now deduce the following proposition.

\begin{prop}\label{prop add part1}
For every $\kappa,\nu\vdash m\geq1,$ and $l\geq1$ such that $\rk(\kappa)=\rk(\nu)+m-l$, we have
\begin{equation}
\mc^\kappa_{\nu,l}=\mc^{\kappa\cup1^n}_{\nu\cup{1^n},l+n}, \hspace{0.5 cm} \text{for every  }n\geq 1.    
\end{equation}
\begin{proof}
First, notice that from \cref{lem bounddeg} we have that $\mc^\kappa_{\nu,l}=\mc^\kappa_{\nu,l}(0)$. On the other hand,
\begin{align*}
\mc^{\kappa\cup1^n}_{\nu\cup{1^n},l+n}(0)&=\sum_{\underset{\ell(\theta)=l+n}{\theta\vdash n+m}}c^{\kappa\cup1^n}_{\nu\cup{1^n},\theta}(0)
\end{align*}
From \cref{lem add part1} item (\textit{2}), the partitions $\theta$ which contribute to this sum are of the form $\theta=\widetilde \theta\cup{1^n}$, where $\widetilde \theta\vdash m$ and $\ell(\widetilde\theta)=l.$ We conclude by using \Cref{lem add part1} item (\textit{1}).
\end{proof}
\end{prop}

\begin{defi}\label{def top coefficients}
Let $\rho$ and $\pi$ be two partitions  of size $r\geq 1$. We define the top coefficient\footnote{This terminology will be justified in \Cref{sec FHA}.} $\tc^\rho_{\pi}$   by $$\tc^\rho_{\pi}:=\mc^{\kappa}_{\nu,l},$$ where   $\kappa $ and $\nu$  are two partitions of the same size $n\geq r+\ell(\pi)$, such that $\kappa=\rho\oplus 1^{n-r}$ (or equivalently $\rho=\kappa-\mathbf{1}$),  $\nu=\pi\cup 1^{n-r}$ and $l$ is such that $n-l+\rk(\pi)=r$. 
 \end{defi}
Note that given \Cref{prop add part1} this definition does not depend on $n$. We consider the matrix of top coefficients defined for any $r\geq 1$ by $\mathcal{Q}^{(r)}:=(\tc^{\rho}_{\pi})_{\pi,\rho\vdash r}$.

\begin{exe}
For $r=3$, the matrix $\mathcal{Q}^{(r)}$ is given by 

\begin{center}
 \begin{tabular}{|c|c|c|c|}
 \hline
 $\pi$ $\backslash$ $\rho$ &  $[3]$ & $[2,1]$ & $[1^3]$\\ \hline
 $[3]$& 4 & 1 & 0  \\ \hline
 $[2,1]$ & 6 & 4 &3  \\ \hline
 $[1^3]$ & 1 & 1 & 1 \\ \hline
   \end{tabular}.
   
\end{center}

\end{exe}

The following theorem will be proved in \Cref{sec FHA} (see also \Cref{thm g}).

\begin{thm}\label{thm Qinv}
The matrix $\mathcal{Q}^{(r)}=(\tc^\rho_\pi)_{\rho,\pi\vdash r}$ is invertible in $\mathbb Z$ for every $r\geq 1.$
\end{thm}

There exists an explicit expression of the top coefficients $\tc^\rho_\pi$, see \cite{BG92,GS98}. However, for the proof of \Cref{thm Qinv}, it will be more natural to consider the algebraic definition of these coefficients and see the matrix $\mathcal{Q}^{(r)}$ as a change-of-basis matrix (see \Cref{prop Q tr mat}).

\section{Proof of \texorpdfstring{\Cref{thm integrality}}{}}\label{sec proof main result}
 The main purpose of this section is to prove that \Cref{thm Qinv} implies \Cref{thm integrality}.
We start by proving the multiplicativity property satisfied by the coefficients $c^\lambda_{\mu,\nu}$ (see also \Cref{eq multiplicativity}): 
\begin{prop}\label{prop multiplicativity}
For every $\lambda,\mu,\nu,\rho \vdash n\geq1$, we have 
    
\begin{equation*}
\sum_{\kappa\vdash n} c^\lambda_{\mu,\kappa}c^\kappa_{\nu,\rho}=\sum_{\theta\vdash n}c^\lambda_{\theta,\rho}c^\theta_{\mu,\nu}.
\end{equation*}
In particular, for every $\lambda,\mu,\nu,\vdash n\geq1$ and $l\geq1$,
\begin{equation*}
\sum_{\kappa\vdash n} c^\lambda_{\mu,\kappa}\mc^\kappa_{\nu,l}=\sum_{\theta\vdash n}\mc^\lambda_{\theta,l}c^\theta_{\mu,\nu}.
\end{equation*}
\begin{proof}
Following \cite{BD21}, we introduce a four parameter version of the coefficients $c^\lambda_{\mu,\nu}$ (see also \cref{ssec k parameters}), defined by
$$\sum_{n\geq0}t^n\sum_{\theta\vdash n}\frac{1}{j^{(\alpha)}_\theta}J^{(\alpha)}_\theta(\mathbf{p})J^{(\alpha)}_\theta(\mathbf{q})J^{(\alpha)}_\theta(\mathbf{r})J^{(\alpha)}_\theta(\mathbf{s})=1+\sum_{n\geq 1}t^n\sum_{\lambda,\mu,\nu,\rho\vdash n}\frac{c^\lambda_{\mu,\nu,\rho}(b)}{z_\lambda(1+b)^{\ell(\lambda)}}p_\lambda q_\mu r_\nu s_\rho,$$ 
where $\mathbf{p},\mathbf{q},\mathbf{r}$ and $\mathbf{s}$ are sequences of power-sum variables as in \cref{ssec coef c}.

The following relation is a consequence of the orthogonality of Jack polynomials (see \cite[Proposition 6.1]{BD21}):
 \begin{equation}\label{eq c four parameters}
   c^\lambda_{\mu,\nu,\rho}=\sum_{\kappa\vdash n}c^\lambda_{\mu,\kappa}c^\kappa_{\nu,\rho}.  
 \end{equation}
On the other hand, it is easy to see from the definition, that the coefficients $c^\lambda_{\mu,\nu,\rho}$ are symmetric in the parameters $\mu$, $\nu$ and $\rho$. In particular, we have $$c^\lambda_{\mu,\nu,\rho}=c^\lambda_{\rho,\mu,\nu}.$$
We conclude by taking the sum over partitions $\rho$ of length $l$.
\end{proof}
\end{prop}

We fix $n>0$. 
 For $0\leq r< n$, we introduce the assertion $\mathcal A_n^{(r)}$:
\begin{align*}\label{eq assertion}
    \mathcal A_n^{(r)}:  &\text{ for every }\lambda,\mu,\kappa\vdash n \text{ such that } \rk(\kappa)=r, \text{ the coefficient }c^\lambda_{\mu,\kappa}\text{ is an integer polynomial in } b.
\end{align*}
Our purpose is to prove prove $\mathcal A_n^{(r)}$ by induction on $r$. We start by the following lemma.
\begin{lem}\label{lem linearsystem}
We fix $1\leq r<n$, and  we assume that the assertions $\mathcal A_n^{(i)}$ hold for $i<r$. 
Let $\lambda,\mu\vdash n$ and let $(\nu,l)$ be a pair satisfying the condition
\begin{equation}\label{condition}
    \nu\vdash n, \text{ }\rk(\nu)<r, \text{   and  } n-l+\rk(\nu)=r\tag{C1}.
\end{equation} 
Then the quantity 
\begin{equation}\label{eq system}
 P^{(r)}_{\lambda,\mu,\nu,l}:=\sum_{\rk(\kappa)=r} c^\lambda_{\mu,\kappa}\mc^\kappa_{\nu,l},  
 \end{equation}
is an integer polynomial in $b$.
\begin{proof}
Note that with the conditions of the proposition, the right hand-side $\sum_{\theta\vdash n} \mc^\lambda_{\theta,l}c^\theta_{\mu,\nu}$ in \Cref{eq multiplicativity}  is an integer polynomial (we use the induction hypothesis and \Cref{thm mar}). This implies that the left hand-side $\sum_{\kappa\vdash n}c^\lambda_{\mu,\kappa}\mc^\kappa_{\nu,l}$ in \Cref{eq multiplicativity}  is an integer polynomial. We conclude using the two following facts: 
\begin{itemize}
    \item if $\rk(\kappa)>r$ then $\mc^\kappa_{\nu,l}=0$.
    Indeed, from \Cref{lem bounddeg}, we get that $$\deg_b(\mc^\kappa_{\nu,l})\leq n-l+\rk(\nu)-\rk(\kappa)=r-\rk(\kappa)<0.$$
    \item if $\rk(\kappa)<r$ then  we know that $c^\lambda_{\mu,\kappa}$ is an integer polynomial from the induction hypothesis.\qedhere
    \end{itemize}
\end{proof}
\end{lem}

For fixed $\lambda,\mu\vdash n$ and  $r<n$,  we get from the previous lemma more equations of type $\eqref{eq system}$  than variables $c^\lambda_{\mu,\kappa}$, where $\rk(\kappa)=r$. In order to obtain a square system, we consider the equations $\eqref{eq system}$ that are indexed by pairs $(\nu,l)$ satisfying the following condition refining \eqref{condition}: \begin{equation}\label{condition 2}
     (\nu,l)=(\pi\cup1^{n-r},l), \text{ where } \pi\vdash r \text{ and } n-l+\rk(\pi)=r.\tag{C2}
 \end{equation}

We denote by $\mathcal{S}^{(r)}_{\lambda,\mu}$ the linear system obtained by taking the equations $\eqref{eq system}$ for $(\nu,l)$ satisfying condition \eqref{condition 2}, and we denote by $\mathcal{Q}^{(r)}_n$ the matrix associated to this system. In other terms $\mathcal{Q}^{(r)}_n=(\mc^\kappa_{\nu,l})$ where indices $\kappa$ of columns are partitions of $n$ of rank $r$, and indices of rows are pairs $(\nu,l)$ satisfying condition \eqref{condition 2}.
Note that this matrix is independent from $\lambda$ and $\mu$. 
In the case $2r\leq n$, one can see that the system $\mathcal{S}^{(r)}_{\lambda,\mu}$ is a square system whose number of columns and number of rows are both equal to the number of partitions of size $r$. Moreover, in this case $\mathcal Q_n^{(r)}$ is the matrix $\mathcal Q^{(r)}$ defined in \cref{sec top coef}. The situation in the case $2r>n$ is more intricate, since $\mathcal{Q}_n^{(r)}$ has less columns. In general, we have the following proposition:

\begin{prop}\label{prop Qn}
For every $1\leq r<n$, the matrix $\mathcal{Q}^{(r)}_n$ is a submatrix of $\mathcal{Q}^{(r)}$ defined in \Cref{sec top coef}, obtained by erasing only columns. More precisely,  $\mathcal{Q}^{(r)}_n=(t^\rho_\pi)$ where the rows index $\pi$ is a partition of $r$, and the columns index $\rho$ is a partition of $r$, such that $\ell(\rho)\leq n-r$.
\begin{proof}
We have the following bijection
\begin{align}\label{bijection}
 \left\{\kappa \text{ partition of $n$ such that } \rk(\kappa) =r\right\} &\xrightarrow[]{\sim} \left\{\rho \text{ partition of $r$ such that } \ell(\rho)\leq n-r\right\}\\
  \kappa&\longmapsto \kappa-\mathbf{1}\nonumber\\
  \rho\oplus1^{n-r}&\text{\reflectbox{$\longmapsto$}}\rho. \nonumber 
\end{align}
Moreover, we recall that $\tc^\rho_\pi=\mc^{\kappa}_{\pi\cup1^{n-r},n-r+\rk(\pi)}$, where $\kappa=\rho\oplus 1^{n-r}$ (see \Cref{def top coefficients}). This concludes the proof.
\end{proof}
\end{prop}

We now prove \Cref{thm integrality}.

\begin{proof}[Proof of \cref{thm integrality}]We prove $\mathcal{A}_{n}^{(r)}$ by induction on $r$.
 For $r=0$, the only partition of rank $0$ is $\kappa=[1^n]$ and we know that $c^\lambda_{\mu,[1^n]}=\delta_{\lambda,\mu}$  for all partitions $\lambda,\mu$, where $\delta_{\lambda,\mu}$ is the Kronecker delta (see \cite[Lemma 3.3]{GJ96b}).
 
  Now we fix $r>0$ and we assume that $\mathcal{A}_{n}^{(j)}$ holds for each $j\leq r-1$. We fix two partitions $\lambda,\mu\vdash n\geq 1$, and we consider the system $\mathcal{S}^{(r)}_{\lambda,\mu}$. It can be written as follows
$$\mathcal{Q}^{(r)}_n X^{(r)}_{\lambda,\mu}=Y_{\lambda,\mu}^{(r)},$$ 
where $Y_{\lambda,\mu}$ is the column vector containing the polynomials $P^{(r)}_{\lambda,\mu,\nu,l}$ for $(\nu,l)$ satisfying condition \eqref{condition 2},
and $X^{(r)}_{\lambda,\mu}$ is the column vector containing $c^\lambda_{\mu,\kappa}$ for $\kappa\vdash n$ of rank $r$. 
We define the column vector $\widetilde X^{(r)}_{\lambda,\mu}:=(x^\lambda_{\mu,\rho})$ for $\rho\vdash r$, where 
$$x^\lambda_{\mu,\rho}:=\left\{ \begin{array}{cc}
    c^\lambda_{\mu,\rho\oplus1^{n-r}} &\text{ if } \ell(\rho)\leq n-r  \\
    0 & \text{otherwise.}
\end{array}\right.$$
Using \cref{prop Qn}, the previous system can be rewritten as follows
$$\mathcal{Q}^{(r)}\widetilde X^{(r)}_{\lambda,\mu}=Y^{(r)}_{\lambda,\mu}.$$
But we know from \Cref{thm Qinv} that $\mathcal{Q}^{(r)}$ is invertible in $\mathbb{Z},$ and since the entries of $Y_{\lambda,\mu}$ are integer polynomials in $b$ (see \Cref{lem linearsystem}), we deduce that this is also the case for the entries of $\widetilde X^{(r)}_{\lambda,\mu}$. Thus the coefficients $c^\lambda_{\mu,\kappa}$ are integer polynomials in $b$, when the partition $\kappa$ has rank $r$. This gives the assertion  $\mathcal{A}_{n}^{(r)}$.
\end{proof}

Note that the previous proof implies that \Cref{eq multiplicativity} allows to recover the coefficients $c^\lambda_{\mu,\nu}$ from their marginal sums. More precisely, we have the following proposition:
\begin{prop}\label{rmq eqsystem}
Fix a real $b$ and let $(y^\lambda_{\mu,\nu})_{\lambda,\mu,\nu\vdash n}$ be a family of numbers indexed by partitions of size $n\geq1$,  satisfying 
$$\left\{
\begin{array}{l}
    y^\lambda_{\mu,[1^n]}=\delta_{\lambda,\mu}  \\
    \sum_{\kappa\vdash n}y^\lambda_{\mu,\kappa}\mc^\kappa_{\nu,l}(b)=\sum_{\theta\vdash n} \mc^\lambda_{\theta,l}(b) y^\theta_{\mu,\nu}.
\end{array}
\right.$$
Then we have that $y^\lambda_{\mu,\nu}=c^\lambda_{\mu,\nu}(b)$ for all partitions $\lambda,\mu,\nu\vdash n$.
\end{prop}

\section{Graded Farahat-Higman Algebra}\label{sec FHA}

This section is dedicated to the proof of \Cref{thm Qinv}.
We start by some notation related to permutations.
If $\sigma$ is a permutation of cyclic type $\lambda$, we define its \textit{reduced cyclic type} as the partition $\lambda-\mathbf{1}$, see \cite{MN13}. Hence if $\sigma\in \mathfrak S_n\subset \mathfrak S_{n+1}...$, the reduced cyclic type of $\sigma$ does not depend on $n$. 
For any partition $\lambda$, we define $C_\lambda(n)\in \mathbb{Z}\mathfrak S_n$ as the sum of all permutations in $\mathfrak S_n$ of reduced cyclic type $\lambda$. Note that $C_\lambda(n)=0$ if $|\lambda|+\ell(\lambda)> n$.

For every $n\geq0$, the family $(C_\lambda(n))_{|\lambda|+\ell(\lambda)\leq n}$ form a basis of the center of the group algebra of $\mathfrak S_n$.
The multiplication in this algebra is given by 
\begin{equation}\label{eq strcuture coef}
    C_\lambda(n) C_\mu(n)=\sum_{|\kappa|+\ell(\kappa)\leq n } \rho^\kappa_{\lambda,\mu}(n) C_\kappa(n),
\end{equation} 
 for some structure coefficients $\rho^\kappa_{\lambda,\mu}(n)$. The latter are linked to the coefficients $c^\lambda_{\mu,\nu}$ evaluated at $b=0$ as follows (see \cite[Proposition 3.1]{GJ96b}): 
 \begin{equation}\label{eq structure coefficients}
  \rho^\kappa_{\lambda,\mu}(n)=\left\{\begin{array}{cc}
     c^{\kappa\oplus 1^{n-|\kappa|}}_{\lambda\oplus 1^{n-|\lambda|},\mu\oplus 1^{n-|\mu|}}(0)  & \text{if } \max\left(|\lambda|+\ell(\lambda),|\mu|+\ell(\mu),|\kappa|+\ell(\kappa)\right)\leq n ,\\
      0 & \text{otherwise.} 
  \end{array}
  \right.   
 \end{equation}
The following proposition is due to Farahat and Higman.
\begin{prop}[\cite{FH59}]\label{prop FH}
The structure coefficients $\rho^\kappa_{\lambda,\mu}$ satisfy the following properties:
\begin{enumerate}
    \item $\rho^\kappa_{\lambda,\mu}=0$ if $|\kappa|>|\lambda|+|\mu|.$
    \item $\rho^\kappa_{\lambda,\mu}$ is independent from $n$ if $|\kappa|=|\lambda|+|\mu|.$
    \item $\rho^\kappa_{\lambda,\mu}$ is a polynomial in $n$ if $|\kappa|\leq|\lambda|+|\mu|.$
\end{enumerate}
\end{prop}

We now prove \cref{lem add part1}.

\begin{proof}[Proof of \cref{lem add part1}]
Item (\textit{1}) can be obtained from \Cref{eq structure coefficients} and \cref{prop FH} item (\textit{2}). 
Item \textit{(2)} is a consequence of \cite[Lemma 3.5]{FH59}.
\end{proof}

We are here interested in the structure coefficients $\rho^{\kappa}_{\lambda,\mu}$ in the case $|\kappa|=|\lambda|+|\mu|$. They are called \textit{the top connection coefficients} of the Farahat-Higman algebra.
To study these coefficients, we consider the graded algebra $\Z_n$ associated to $Z(\mathbb Z\mathfrak S_n)$ with respect to the filtration $$\deg(C_\lambda(n))=|\lambda|.$$
We denote by $\mathfrak c_\lambda(n)$ the image of $C_\lambda(n)$ in $\Z_n.$
Concretely,
$\Z_n=\bigoplus_{1\leq r\leq n-1}\Z^{(r)}_n,$
where $$\Z_n^{(r)}:=\Span_\mathbb{Z}\left\{\cf_\lambda(n);\lambda \vdash r\text{ and }\ell(\lambda)\leq n-r\right\},$$ and  the multiplication in $\Z_n$ is defined by
$$\cf_\lambda(n) \cf_\mu(n)=\sum_{\kappa\vdash |\lambda|+|\mu|}\rho^{\kappa}_{\lambda,\mu}(n)\cf_\kappa(n).$$
Note that compared to \Cref{eq strcuture coef}, we keep only the top degree terms.
The graded algebra $\Z_n$ comes with a linear isomorphism $\phi_n: Z(\mathbb Z\mathfrak S_n)\xrightarrow[]{\sim}\Z_n,$ that sends $C_\lambda(n)$ to $\cf_\lambda(n)$ (which is obviously not an algebra isomorphism).

Since the structure coefficients in $\Z_n$ are independent from $n$, we can define a family of $\mathbb Z$-algebra morphisms:
\begin{align*}
\psi_n:\Z_{n+1}\longrightarrow &\Z_n \\
\cf_\lambda(n+1)\longmapsto &\left\{
\begin{array}{cc}
     \cf_\lambda(n)& \text{ if }|\lambda|+\ell(\lambda)\leq n,  \\
     0& \text{otherwise.} 
\end{array}
\right.
\end{align*}

Let $\Z_\infty:=\varprojlim \Z_n$ be the projective limit of the $\Z_n$'s, and let the limit $\cf _\lambda:=\varprojlim \cf_\lambda(n)$ in $\Z_\infty$.
We also define $\mathcal{Z}_\infty^{(r)}:=\Span _\mathbb Z\{\cf_\lambda;\lambda \vdash r\}$, hence $\Z_\infty=\bigoplus_{r\geq 1}\Z^{(r)}_\infty$ (see \cite[Example 24, page 131]{Mac} for more details about the construction of the algebra $\Z_\infty$).  

We define for every $r \geq 1$, 
$\ff_r(n):=\sum_{\lambda\vdash r}\cf_\lambda(n),$
and for any partition $\mu$,
$\ff_\mu(n):=\prod_i\ff_{\mu_i}(n).$
We also define  
$\mathfrak g_\pi(n):=\cf_{\pi-\mathbf{1}}(n)\ff_{\ell(\pi)}(n)$, for any partition $\pi$.
Note that $$\deg(\ff_\mu(n))=\deg(\gf_\mu(n))=|\mu|.$$ Finally we define the limits $\ff_\mu:=\varprojlim \ff_\mu(n)$ and $\gf_\pi:=\varprojlim \gf_\pi(n)$ in $\Z_\infty.$
We have the following theorem due to Farahat and Higman.
\begin{thm}[\cite{FH59}]\label{thm FH}
For every $r\geq0$, $(\ff_\lambda)_{\lambda\vdash r}$ is a $\mathbb Z$-basis of $\Z^{(r)}$.
\end{thm}

The following lemma relates the top coefficients $\tc^\rho_\pi$ (see \Cref{def top coefficients}) to extraction coefficients in the graded algebra.
\begin{lem}
For all partitions $\rho,\pi\vdash r\geq 0$, we have $\tc^\rho_\pi=[\cf_{\rho}]\gf_{\pi}.$
\begin{proof}
From the definition 
$\tc^\rho_\pi=\mc^\kappa_{\nu,l},$
where   $\kappa $ and $\nu$  are the two partitions of the same size $n\geq r+\ell(\rho)$, such that $\kappa=\rho\oplus1^{n-r}$,  $\nu=\pi\cup1^{n-r}$ and $l$ is such that $n-l+\rk(\nu)=r$. On the other hand, 
$$\mc^\kappa_{\nu,l}=[C_{\kappa-\mathbf{1}}(n)]\left(\sum_{\lambda\vdash n-l}C_\lambda(n)\right)C_{\nu-\mathbf{1}}(n),$$
(we use here \Cref{eq structure coefficients} and the fact that under the condition $n-l+\rk(\nu)=r=\rk(\kappa)$ we have $\mc^\kappa_{\nu,l}=\mc^\kappa_{\nu,l}(0)).$
Since this is a top degree coefficient extraction, we can consider it in the graded algebra $\Z_n$:
$$\tc^\rho_\pi=\mc^\kappa_{\nu,l}=[\cf_\rho(n)]\ff_{n-l}(n)\cf_{\nu-\mathbf{1}}(n).$$
To conclude note that $\nu-\mathbf{1}=\pi-\mathbf{1},$ and that $n-l=r-\rk(\nu)=r-\rk(\pi)=\ell(\pi).$ 
\end{proof}
\end{lem}

We deduce the following proposition.

\begin{prop}\label{prop Q tr mat}
For every $r\geq 0$, the matrix $(\mathcal{Q}^{(r)})^{T}$ is the matrix of $(\gf_\pi)_{\pi \vdash r}$ in the basis $(\cf_\rho)_{\rho \vdash r}.$

\end{prop}

Hence, our goal is to prove the following theorem that implies \Cref{thm Qinv}:
\begin{thm}\label{thm g}
For every $r\geq0$, the family $(\gf_\pi)_{\pi\vdash r}$ is a $\mathbb Z$-basis for $\Z^{(r)}_\infty$.
\end{thm}
The rest of this section is dedicated to the proof of this theorem. Let us explain the steps of this proof;
instead of studying the matrix of $(\gf_\pi)$ in the basis $(\cf_\rho)$, we consider its matrix $\mathcal N$ in the basis $(\ff_\lambda)$. This matrix satisfy some block triangularity property. Its diagonal blocks are submatrices of $\mathcal M$, defined as the matrix of $(\cf_\rho)$ in the basis $(\ff_\lambda)$. In order to prove that these submatrices of $\mathcal M$ are invertible\footnote{Note that \Cref{thm FH} implies that the matrix $\mathcal{M}^{(r)}$ has determinant $\pm 1$, however we will need to obtain this result for some submatrices of $\mathcal{M}^{(r)}$ called South-East blocks, see \cref{cor M}.} in $\mathbb{Z}$, we will need to introduce an intermediate basis $(\mf_\mu)$ (see definition below). This basis has the property that its transition matrices to both bases $(\cf_\rho)$ and $(\ff_\lambda)$ (denoted respectively $\mathcal{L}$ and $\mathcal{U}$) are triangular (\cref{prop M1} and \Cref{Thm M2}).

The following diagram illustrates the different basis of $\mathcal{Z}^{(r)}_\infty$ that will be useful for the proof of \Cref{thm g}, and the associated transition matrices.
\begin{center}
\begin{tikzcd}    
 
 (\cf_\rho)_{\rho\vdash r}\arrow[dd,swap,"\mathcal{M}^{(r)}"] \ar[ddrr,bend right=10,"\mathcal{L}^{(r)}" near end]  &&(\gf_\pi)_{\pi\vdash r}\arrow[ll,swap,"{(Q^{(r)})^T}"]  \arrow[lldd,"\mathcal{N}^{(r)}" near start,bend right=25]\\
 \\
   (\ff_\lambda)_{\lambda\vdash r}    &&(\mf_\mu)_{\mu\vdash r}\arrow{ll}{\mathcal{U}^{(r)}}\\
 \end{tikzcd}
 \end{center}

We denote by $(\mathcal{J}_i)_{i\geq2}$ the Jucys-Murphy elements:
$$\mathcal{J}_i:=(1,i)+....+(i-1,i)\in\mathbb{Z}\mathfrak{S}_n, \text{ for every }n\geq i.$$
We denote for any symmetric function $f$, 
$$f(\Xi_n)=f(\mathcal{J}_2,\mathcal{J}_3,...,\mathcal{J}_n,0,0,..).$$

The evaluation of the elementary symmetric function in the Jucys-Murphy elements has the following expression (see \cite{J74}):
$$e_l(\Xi_n)=\sum_{\lambda\vdash l}C_\lambda(n), \text{ for every }l\geq 1.$$ Hence the image of $e_\lambda(\Xi_n)$ in $\Z_n$ is given by $\phi_n\left(e_\lambda(\Xi_n)\right)=\ff_\lambda(n).$ Since the elementary functions form a basis of the symmetric functions algebra, we get that $f(\Xi_n)\in Z(\mathbb Z\mathfrak S_n)$ for any symmetric function $f$.
We now consider the evaluation of the monomial symmetric functions in the Jucys-Murphy elements.

\begin{defi}
For any partition $\mu\vdash r \geq0$, we define 
$\mf_\mu(n)\in\Z_n$ by 
$$\mf_\mu(n):= \phi_n\left(m_\mu(\Xi_n)\right),$$ 
where $m_\mu$ denotes the monomial symmetric function associated to $\mu$. We also introduce their limit  in $\Z_\infty$;
$\mf_\mu:=\varprojlim \mf_\mu(n).$
\end{defi}
Note that since the elementary functions form a basis of the symmetric functions algebra, $\mf_\mu(n)$
is a linear combination of $\ff_\lambda(n)$, and the previous limit is well defined. The elements $\mf_\mu(n)$ have been studied in \cite{MN13}, and an explicit expression of their developments in the basis $(\cf_\rho(n))$ has been given.

In the following, we will study some triangularity properties of the different transition matrices in $\Z_\infty$, where the three families $(\cf_\rho)$, $(\mf_\mu)$ and $(\gf_\pi)$ will be indexed with partitions with respect to the total order $\leq$, and the basis $(\ff_\lambda)$ will be indexed with partitions with respect to the dual order $\leq'$ (see \Cref{sec partitions} for the definition of these orders). Moreover, we order rows and columns in increasing order (see examples below).

\begin{prop}\label{prop M1}
The matrix $\mathcal{U}^{(r)}$ of $(\mf_\mu)_{\mu\vdash r,\leq}$ in the basis  $(\ff_\lambda)_{\mu\vdash r,\leq'}$ is upper triangular, with diagonal coefficients equal to 1, i.e
$$[\ff_\lambda]\mf_\mu=0 \text{ if } \mu<\lambda' \text{ and } [\ff_{\mu'}]\mf_\mu=1.$$
\begin{exe}
For $r=3$, the matrix $\mathcal{U}^{(r)}$ is given by 

\begin{center}
\begin{tabular}{|c|c|c|c|}
\hline
     $\ff_\lambda \backslash \mf_\mu$ &  $[1^3]$ & $[2,1]$ &[3] \\\hline
     $[3]$  & 1  & -3 &  3 \\\hline
     $[2,1]$ & 0 & 1  &  -3 \\\hline
     $[1^3]$ & 0 & 0 & 1  \\\hline
\end{tabular}.
\end{center}
\end{exe}
 \begin{proof}
 We know that the matrix of the monomial functions $m_\mu$ in the elementary symmetric functions basis  $e_\lambda$ has the triangularity property (see \cite[Equation (2.3)]{Mac}):
 $$m_\mu=\sum_{\lambda\vdash r}u_{\lambda,\mu}e_\lambda,$$ 
 where $u_{\lambda,\mu}$ are integers such that $u_{\mu',\mu}=1$ and $u_{\lambda,\mu}=0$ if $\mu<\lambda'.$ 

 By evaluating at the Jucys-Murphy elements and applying $\phi_n$ we obtain 
 $$\mf_\mu(n)=\sum_{\mu\vdash r}u_{\lambda,\mu}\ff_\lambda(n).$$
 We conclude by taking the limit in $\Z_\infty$.
 \end{proof}
\end{prop}

The following theorem is due to Matsumoto and Novak.
\begin{thm}\cite[Theorem 2.4]{MN13}\label{Thm M2}
For every $r\geq0$, the matrix $\mathcal{L}^{(r)}$  of  $(\cf_\rho)_{\rho\vdash r,\leq}$ in $(\mf_\mu)_{\mu\vdash r,\leq}$
is lower triangular, with diagonal coefficients equal to 1, i.e
$$[\mf_\mu]\cf_\rho=0 \text{ if } \mu<\rho \text{ and } [\mf_\rho]\cf_\rho=1.$$
\end{thm}
\begin{exe}
For $r=3$, $\mathcal{L}^{(r)}$ is given by

\begin{center}
\begin{tabular}{|c|c|c|c|}\hline
     $\mf_\mu \backslash \cf_\rho$ &  $[1^3]$ & $[2,1]$ &[3] \\\hline
     $[3]$  & 1  & 0 & 0  \\\hline
     $[2,1]$ & -1 & 1  &  0 \\\hline
     $[1^3]$ & 2 & -3 & 1  \\\hline
\end{tabular}.
\end{center}

\end{exe}

For every $r\geq0$, we define $\mathcal{M}^{(r)}$ as the matrix of $(\cf_\rho)_{\rho\vdash r,\leq}$ in $(\ff_\lambda)_{\lambda\vdash r,\leq'}$. Hence 

$$\mathcal{M}^{(r)}=\mathcal{U}^{(r)}\mathcal{L}^{(r)}.$$

For every $1\leq i\leq r$, we define $\mathcal{M}^{(r,i)}$ as the submatrix of $\mathcal{M}^{(r)}$ obtained by keeping the indices $(\lambda,\rho)$ such that the rows index $\lambda$ satisfies $\lambda_1\leq i$ and the columns index $\rho$ satisfies $\ell(\rho)\leq i$. This matrix is a South-East block of  $\mathcal{M}^{(r)}$ (this is a consequence of item (2) in \cref{def orders}).
\begin{exe}\label{Ex M}
For $r=3$, 

$\mathcal{M}^{(3)}$ is given by
\begin{tabular}{|c|c|c|c|}\hline
     $\ff_\lambda \backslash \cf_\rho$& $[1^3]$ & $[2,1]$ &[3] \\\hline
     $[3]$  & 10  & -12 &  3 \\\hline
     $[2,1]$ & -7 & 10 &  -3 \\\hline
     $[1^3]$ & 2 & -3 & 1  \\\hline
\end{tabular}
and $\mathcal{M}^{(3,2)}$ by
\begin{tabular}{|c|c|c|}\hline
     $\ff_\lambda\backslash \cf_\rho$ & $[2,1]$ &[3] \\\hline
     $[2,1]$  & 10 &  -3 \\\hline
     $[1^3]$  & -3 & 1  \\\hline
\end{tabular}.

\end{exe}

Similarly, we define $\mathcal{U}^{(r,i)}$ as the submatrix of $\mathcal{U}^{(r)}$ obtained by keeping the indices $(\lambda,\mu)$ such that $\lambda_1\leq i$ and $\ell(\mu)\leq i$, and $\mathcal{L}^{(r,i)}$ as the submatrix of $\mathcal{L}^{(r)}$ obtained by keeping the indices $(\mu,\rho)$ such that $\ell(\mu)\leq i$ and $\ell(\rho)\leq i$.
\begin{cor}\label{cor M}
For every $r\geq1$ and $1\leq i\leq r$, the submatrix  $\mathcal{M}^{(r,i)}$ has determinant 1.
\begin{proof}
Using the triangularity properties of  $\mathcal{L}^{(r)}$ and $\mathcal{U}^{(r)}$ one can check that 

$$\mathcal{M}^{(r,i)}=\mathcal{U}^{(r,i)}\mathcal{L}^{(r,i)}, $$
and that $\det\left(\mathcal{U}^{(r,i)}\right)=\det\left(\mathcal{L}^{(r,i)}\right)=1$.
\end{proof}
\end{cor}

We now prove \cref{thm g}.

\begin{proof}[Proof of \cref{thm g}]
Let $\mathcal N^{(r)}$ be the matrix of $(\gf_\pi)_{\pi\vdash r,\leq}$ in $(\ff_\lambda)_{\lambda\vdash r,\leq'}$. Our purpose is to prove that $\det(\mathcal{N}^{(r)})=1$.
Let  $\mathcal N^{(r)}_{\lambda,\pi}$ denote the coefficients of $\mathcal N^{(r)}$, namely
\begin{equation}\label{eq g in f}
\gf_\pi=\sum_{\lambda\vdash r}\mathcal N^{(r)}_{\lambda,\pi}\ff_\lambda.    
\end{equation}
We start by proving that $\mathcal N^{(r)}$ is block-upper triangular. More precisely, if  $\mathcal N^{(r)}_{\lambda,\pi}\neq0$ then $\lambda_1\geq\ell(\pi).$ Indeed,
from the definition of $\gf_\pi$ we have
$\gf_\pi=\cf_{\pi-\mathbf{1}}\ff_{\ell(\pi)}.$
We write $\cf_{\pi-\mathbf{1}}$ in the basis $(\ff_\kappa)_{\kappa\vdash r-\ell(\pi)}$,
$$\cf_{\pi-\mathbf{1}}=\sum_{\kappa\vdash r-\ell(\pi)}\mathcal{M}^{(r-\ell(\pi))}_{\kappa,\pi-\mathbf{1}}\ff_\kappa\text{,}$$
where $(\mathcal{M}^{(r-\ell(\pi))}_{\kappa,\nu})_{\kappa,\nu\vdash r-\ell(\pi)}$ are the coefficients of the matrix $\mathcal{M}^{(r-\ell(\pi))}.$
Multiplying the last equation by $\ff_{\ell(\pi)}$ and comparing it to \Cref{eq g in f}, we get  
\begin{equation}\label{eq af}
  \mathcal{N}^{(r)}_{\lambda,\pi}=\mathcal{M}^{(r-\ell(\pi))}_{\lambda\backslash {\ell(\pi)},\pi-\mathbf{1}}.  
\end{equation}

In particular if  $\mathcal{N}^{(r)}_{\lambda,\pi}\neq0$ then  $\ell(\pi) \text{ is a part of } \lambda$, and then necessarily $\ell(\pi)\leq \lambda_1$.

This proves that that $\mathcal N^{(r)}$ is a block-upper triangular matrix. We define for $1\leq i\leq r$  the diagonal block $\mathcal N^{(r,i)}$ as the submatrix that contains the coefficients $\mathcal{N}^{(r)}_{\lambda,\pi}$ for $\ell(\pi)=\lambda_1=i$. Moreover, we have a bijection between pairs of partitions $(\lambda,\pi)$ of size $r$ such that $\ell(\pi)=\lambda_1=i$ and pairs of partitions $(\nu,\kappa)$ of size $r-i$ and such that $\nu_1\leq i$ and $\ell(\kappa)\leq i$, given by:
$$(\lambda,\pi)\longmapsto (\lambda\backslash\lambda_1,\pi-\mathbf{1})$$
Using \Cref{eq af}, we deduce that 
$$\mathcal{N}^{(r,i)}=\mathcal{M}^{(r-i,i)}.$$
This implies that $\det(\mathcal{N}^{(r,i)})=1$. Since $\mathcal{N}^{(r)}$ is block triangular, we deduce that $\det\left(\mathcal{N}^{(r)}\right)=1$, and that $(\gf_\pi)_{\pi\vdash r}$ is a basis for $\Z_\infty^{(r)}.$
\end{proof}

\begin{exe}
We give here the matrix $\mathcal{N}^{(r)}$ for $r=5$. The diagonal blocks, defined by $\ell(\pi)=\lambda_1$, are colored in gray. Note that for $i=2$, $\mathcal{N}^{(5,2)}$ is equal to $\mathcal{M}^{(3,2)}$ given in  \Cref{Ex M}.
\definecolor{mygray}{gray}{0.8}
\begin{center}
\begin{tabular}{|c|c|c|c|c|c|c|c|}
\hline
 $\ff_\lambda$ $\backslash$ $\gf_\pi$ & $[1^5]$  & $[2,1^3]$ & $[2^2,1]$ & $[3,1^2]$ &  $[3,2]$ &  $[4,1]$   & $[5]$\\ \hline
     $[5]$    & \cellcolor{mygray}1& 0 & 0& 0& 0 &  0 & 0 \\ \hline
     $[4,1]$  &  & \cellcolor{mygray}1 & 0& 0&  0& 0 & -4\\ \hline
     $[3,2]$  &  &   &\cellcolor{mygray} 3 &\cellcolor{mygray}-1&-12& 3 &  0\\ \hline
     $[3,1^2]$&  &   &\cellcolor{mygray} -2&\cellcolor{mygray} 1&  0&  0 &  4 \\ \hline
     $[2^2,1]$&  &   &  &  &\cellcolor{mygray} 10&\cellcolor{mygray}-3& 2 \\ \hline
     $[2,1^3]$&  &   &  &  & \cellcolor{mygray}-3&\cellcolor{mygray}1 & -4\\ \hline
     $[1^5]$  &  &   &  &  &   &  & \cellcolor{mygray} 1\\ \hline
       \end{tabular}
    
\end{center}

\end{exe}

\section{The cases \texorpdfstring{$b=0$}{} and \texorpdfstring{$b=1$}{} in the Matching-Jack conjecture via the multiplicativity property}\label{sec b=0 and b=1}
In this section, we give a combinatorial interpretation of the multiplicativity property of \Cref{eq multiplicativity} in terms of matchings in the special cases $b=0$ and $b=1$. We use this interpretation to give a new proof of the Matching-Jack conjecture in these two cases. Unlike the classical proof given in \cite{GJ96a}, the approach that we use here does not use representation theory. A key step in the proof is a combinatorial  version of \cref{thm mar} that we formulate here in the cases $b=0$ and $b=1$ (see \cref{thm marg b=0 and b=1}). 

We start by introducing matchings, which are the central combinatorial objects in this section, and we use them to give the combinatorial formulation of the Matching-Jack conjecture.
\subsection{Matchings}\label{ssec matchings}

For $n\geq1$, we will consider matchings on the set $\mathcal{A}_n:=\{1,\widehat{1},...,n,\widehat{n}\}$, that is a set partition of $\mathcal{A}_n$ into pairs. A matching $\delta$ on $\mathcal{A}_n$ is \textit{bipartite} if each one of its pairs is of the form $\{i,\widehat j\}$.
Given two matchings $\delta_1$ and $\delta_2$, we consider the partition $\Lambda(\delta_1,\delta_2)$ of $n$ obtained by reordering the half-sizes of the connected components of the graph formed by $\delta_1$ and $\delta_2$.
\begin{exe}
We consider the two matchings $\delta_1=\left\{\{1,\widehat 3\},\{3,\widehat4\},\{4,\widehat 1\},\{2,\widehat 2\}\right\}$ and  $\delta_2=\left\{\{1,3\},\{\widehat 3,\widehat4\},\{4,\widehat 2\},\{2,\widehat 1\}\right\}$ on the set $\mathcal{A}_4$. Then $\delta_1$ is bipartite and $\delta_2$ is not. Moreover, one has $\Lambda(\delta_1,\delta_2)=[2,2].$ 
 
\end{exe}

The Matching-Jack conjecture has the following combinatorial version (see \cite[Conjecture 4.2]{GJ96b}).
\begin{conj_MJ2}[\cite{GJ96b}]
Fix a partition $\lambda\vdash n\geq 1$, and two bipartite matchings $\delta_1$ and $\delta_2$ on $\mathcal{A}_n$ such that $\Lambda(\delta_1,\delta_2)=\lambda$. There exists a statistic $\vartheta_\lambda$ on the matchings of $\mathcal{A}_n$ with non-negative integer values, such that
\begin{itemize}
    \item $\vartheta_\lambda(\delta)=0$ if and only if $\delta$  is bipartite.
    \item for any partitions $\mu,\nu\vdash n$, we have 
$$c^\lambda_{\mu,\nu}=\sum_{\delta}b^{\vartheta_\lambda(\delta)},$$
where the sum is taken over matchings $\delta$ of $\mathcal{A}_n$, such that $\Lambda(\delta_1,\delta)=\mu$ and $\Lambda(\delta,\delta_2)=\nu$.

\end{itemize}
\end{conj_MJ2}
Given the results for $b=0$ and $b=1$ (see \cref{thm b=0 and b=1} below), it is easy to see that this formulation of the Matching-Jack conjecture is equivalent to the one given in the introduction.

\subsection{Multiplicativity for Matchings}

Let $\lambda\vdash n\geq1$ and let $\delta_1$ and $\delta_2$ two bipartite matchings on $\mathcal{A}_n$ such that $\Lambda(\delta_1,\delta_2)=\lambda.$ We define for any partitions $\mu,\nu\vdash n$ the two quantities
$$a^\lambda_{\mu,\nu}=\left|\left\{\delta\text{ such that }\Lambda(\delta_1,\delta)=\mu\text{ and }\Lambda(\delta,\delta_2)=\nu\right\}\right|,$$
and 
$$\widetilde a^\lambda_{\mu,\nu}=\left|\left\{\delta\text{ bipartite matching such that }\Lambda(\delta_1,\delta)=\mu\text{ and }\Lambda(\delta,\delta_2)=\nu\right\}\right|.$$

One can see that these quantities do not depend on the bipartite matchings $\delta_1$ and $\delta_2$ (see e.g. \cite{HSS92}).
They satisfy the multiplicativity property:
\begin{prop}\label{prop multuplicativity b=0 and b=1}
For all partitions $\lambda,\mu,\nu,\rho\vdash n\geq1$, we have  
\begin{equation}\label{eq a 1}
\sum_{\kappa\vdash n} a^\lambda_{\mu,\kappa}a^\kappa_{\nu,\rho}=\sum_{\theta\vdash n}a^\lambda_{\theta,\rho}a^\theta_{\mu,\nu}.
\end{equation}
\begin{equation}\label{eq a 2}
\sum_{\kappa\vdash n} \widetilde a^\lambda_{\mu,\kappa}\widetilde a^\kappa_{\nu,\rho}=\sum_{\theta\vdash n}\widetilde a^\lambda_{\theta,\rho}\widetilde a^\theta_{\mu,\nu}.
\end{equation}
\begin{proof}
We fix two bipartite matchings $\delta_1$ and $\delta_2$ of $\mathcal{A}_n$ such that $\Lambda(\delta_1,\delta_2)=\lambda.$
We introduce the set
$$\F^\lambda_{\mu,\nu,\rho}:=\left\{(\delta,\delta')\text{ such that }\Lambda(\delta_1,\delta)=\mu,\Lambda(\delta,\delta')=\nu,\text{ and }\Lambda(\delta',\delta_2)=\rho\right\}.$$

For every $(\delta,\delta')\in\F^\lambda_{\mu,\nu,\rho}$ we have a diagram 
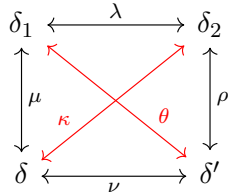
\begin{center}
 \begin{tikzcd}    
  \delta_1 \arrow[rr,"\lambda",leftrightarrow]
  \arrow[dd,"\mu",leftrightarrow]\arrow[rrdd,red,near end,"\theta",leftrightarrow]
  &&\delta_2 \arrow[lldd,red,near end,swap,"\kappa",leftrightarrow] \arrow[dd,"\rho",leftrightarrow]\\
  \\
 \delta && \delta' \arrow[ll,"\nu",leftrightarrow]
 \end{tikzcd}
\end{center}
for some partitions $\theta$ and $\kappa$ of $n$, where the arrow between two matchings $\delta_a$ and $\delta_b$ is labelled by a partition $\xi$ if $\Lambda(\delta_a,\delta_b)=\xi$.
We have two ways to count the number of pairs $(\delta,\delta')$ in $\F^\lambda_{\mu,\nu,\rho}$; either we start by choosing $\delta$ (this corresponds to the left-hand side in \Cref{eq a 1}) or we start by choosing $\delta'$ (this corresponds to the right-hand side in \Cref{eq a 1}).

Similarly, we obtain \Cref{eq a 2} by considering the set 
$$\widetilde \F^\lambda_{\mu,\nu,\rho}:=\left\{(\delta,\delta')\text{ bipartite matchings such that }\Lambda(\delta_1,\delta)=\mu,\Lambda(\delta,\delta')=\nu,\text{ and }\Lambda(\delta',\delta_2)=\rho\right\}.\qedhere$$
\end{proof}
\end{prop}
\subsection{The Matching-Jack conjecture for \texorpdfstring{$b=0$}{} and \texorpdfstring{$b=1$}{}}
The positivity result of \cref{thm mar} on the marginal coefficients $\mc^\lambda_{\mu,l}$ has the following combinatorial version in terms of matchings, that we formulate here in the cases $b=0$ and $b=1$ (this is a special case in \cite[Theorem 1.4]{BD21}).

\begin{thm}[\cite{CD20,BD21}]\label{thm marg b=0 and b=1}
For any integers $n,l\geq 1$ and partitions $\lambda,\mu\vdash n$ we have:
$$\mc^\lambda_{\mu,l}(1)=\sum_{\ell(\nu)=l}a^\lambda_{\mu,\nu}, \text{ and } \mc^\lambda_{\mu,l}(0)=\sum_{\ell(\nu)=l}\widetilde a^\lambda_{\mu,\nu}.$$
\end{thm}
Further details about the proof of this theorem are provided in \cref{rmq proof mar b=0 and b=1}.

In the following, we use the previous theorem and the multiplicativity property to give a new proof for the cases $b=0$ and $b=1$ in the Matching-Jack conjecture.

\begin{thm}\label{thm b=0 and b=1}
For all partitions $\lambda,\mu,\nu\vdash n\geq 1$, we have 
$$c^\lambda_{\mu,\nu}(1)=a^\lambda_{\mu,\nu}, \text{ and } c^\lambda_{\mu,\nu}(0)=\widetilde a^\lambda_{\mu,\nu}.$$
\begin{proof}
Taking the sum over partitions $\rho$ of length $l$ in Equations \eqref{eq a 1} and \eqref{eq a 2} and using \cref{thm marg b=0 and b=1}  we obtain:
\begin{equation*}
\sum_{\kappa\vdash n} a^\lambda_{\mu,\kappa}\mc^\kappa_{\nu,l}(1)=\sum_{\theta\vdash n}\mc^\lambda_{\theta,l}(1)a^\theta_{\mu,\nu}.
\end{equation*}
\begin{equation*}
\sum_{\kappa\vdash n} \widetilde a^\lambda_{\mu,\kappa}\mc^\kappa_{\nu,l}(0)=\sum_{\theta\vdash n} \mc^\lambda_{\theta,l}(0)\widetilde a^\theta_{\mu,\nu}.
\end{equation*}
We use \cref{rmq eqsystem} to conclude.
\end{proof}
\end{thm}

\begin{rmq}
Let $n\geq1$ and let $(\vartheta_\lambda)_{\lambda\vdash n}$ be a family of statistic on matchings of $\mathcal{A}_n$. For all partitions $\lambda,\mu,\nu\vdash n$, we define 
$$y^\lambda_{\mu,\nu}:=\sum_{\delta\in\mathfrak F^\lambda _{\mu,\nu}}b^{\vartheta_\lambda(\delta)}.$$ 
One possible way to prove that the family $(\vartheta_\lambda)$ is a solution for the Matching-Jack conjecture, is to prove that: 
\begin{enumerate}
    \item $y^\lambda_{\mu,\nu}$ satisfies a multiplicativity property of the form of \Cref{eq multiplicativity}.
    \item $y^\lambda_{\mu,\nu}$ is a solution for the conjecture in the case of marginal sums, namely
    $$\sum_{\ell(\nu)=l}y^\lambda_{\mu,\nu}=\mc^\lambda_{\mu,l}.$$
\end{enumerate}
  We can think of the two previous properties as a generalization of \cref{prop multuplicativity b=0 and b=1} and \cref{thm marg b=0 and b=1} given in this section in the case $b\in\{0,1\}$. From the work of \cite{CD20,BD21}, the second property has been verified for some family of statistics $\vartheta$. Unfortunately, these statistics do not seem to satisfy the first property.
\end{rmq}

\begin{rmq}\label{rmq proof mar b=0 and b=1}
The proof of \cref{thm marg b=0 and b=1} is simpler than the general case of \cite[Theorem 1.4]{BD21}. Indeed, as discussed in the introduction, the Matching-Jack conjecture in the marginal case given in \cite{BD21} has been deduced from an analog result for the $b$-conjecture established in \cite{CD20}. The latter result has been obtained by proving that the function $\tau_b$  satisfies a family of partial differential equations. A key step in the proof consists in proving that the differential operators defining these equations satisfy some commutation relations. These relations have a simple combinatorial proof when $b\in\{0,1\}$ (see \cite[Section 4.3]{CD20}). Moreover, in the special cases $b=0$ and $b=1$, the equivalence between the $b$-conjecture and the Matching-Jack conjecture is well understood since we have an  encoding of bipartite maps with matchings  (see e.g. \cite{GJ96a,BD21}). The connection  between the two conjectures established in \cite{BD21} in the  "marginal sums" case is more intricate. 

\end{rmq}

\section{Some consequences of the main result}
\label{sec cumulants}
\subsection{A partial integrality result  for the \texorpdfstring{$b$}{}-conjecture}\label{ssec bconj}
We consider a "connected" version $d^\lambda_{\mu,\nu}$ of the coefficients $c^\lambda_{\mu,\nu}$, defined by 
\begin{equation}\label{eq d}
  \log\left(\tau_b(t,\mathbf{p},\mathbf{q},\mathbf{r})\right)=\sum_{n\geq1}\sum_{\lambda,\mu,\nu\vdash n}t^n \frac{d^\lambda_{\mu,\nu}}{z_\lambda(1+b)^{\ell(\lambda)}}p_\lambda q_\mu r_\nu.  
\end{equation}

These coefficients are related to the coefficients $h^\lambda_{\mu,\nu}$, introduced by Goulden and Jackson in \cite{GJ96b}, by 
\begin{equation}\label{equation h}
h^\lambda_{\mu,\nu}=\frac{n d^\lambda_{\mu,\nu}}{z_\lambda(1+b)^{\ell(\lambda)-1}},    
\end{equation}

for every $\lambda,\mu,\nu\vdash n\geq1$.
The coefficients $h^\lambda_{\mu,\nu}$ are the main object of the $b$-conjecture.
\begin{b-conj}[\cite{GJ96b}]
For all partitions $\lambda,\mu,\nu\vdash n\geq 1$, the coefficient $h^\lambda_{\mu,\nu}$ is a polynomial in $b$ with non-negative integer coefficients.
\end{b-conj}
Do\l{}{\k{e}}ga and Féray have deduced the polynomiality of the coefficients $h^\lambda_{\mu,\nu}$ from the polynomiality of $c^\lambda_{\mu,\nu}$.
\begin{thm}[\cite{DF17}]\label{thm polynomiality h}
For all partitions $\lambda,\mu,\nu\vdash n\geq 1$, the coefficient $h^\lambda_{\mu,\nu}$ is a polynomial in $b$.
\end{thm}
In this section we prove that the coefficients $d^\lambda_{\mu,\nu}$ are integer polynomials in $b$ (see \cref{cor integrality d}). Using \cref{thm polynomiality h}, this implies that $\frac{z_\lambda}{n}h^\lambda_{\mu,\nu}$ is an integer polynomial in $b$, however we do not have information about the divisibility of the coefficients by $\frac{z_\lambda}{n}$.

\begin{defi}
Fix a set $S$. A \textit{set-partition} of $S$ is an unordered family of non-empty disjoint subsets of $S$ whose union is $S$. We denote a by $\mathcal{P}(S)$ the set of set-partitions of $S$.  If $\pi$ is a set-partition of $\lambda$ into $s$ parts then $s$ is called \textit{the size} of $\pi$ and denoted $|\pi|$. For every $s\geq 1$, we denote by $\mathcal{P}_s(S)$ the set of set-partitions of $S$ of size $s$.
\end{defi}

For any integer $n\geq 1$, we set $\llbracket n\rrbracket:=\left\{1,2,...,n\right\}$.
Fix a partition $\lambda$. 
For any subset $B\subset \llbracket\ell(\lambda)\rrbracket$, we denote by $\lambda_B:=[\lambda_i;i \in B]$ the partition obtained by keeping the parts of $\lambda$ with an index in $B$.

The following is a variant of \cite[Lemma 5.2]{DF17} where it is formulated in terms of cumulants. The proof is quite similar.

\begin{lem}\label{lem c-d}

For any partitions $\lambda,\mu,\nu\vdash n\geq1$, we have 
\begin{equation}\label{eq lem c-d}
d^\lambda_{\mu,\nu}=\sum_{\pi\in\mathcal{P}\left(\llbracket\ell(\lambda)\rrbracket\right)}\sum_{(\mu^B)_{B\in\pi},(\nu^B)_{B\in\pi}}\prod_{B\in\pi}(-1)^{|\pi|-1}(|\pi|-1)!c ^{\lambda_B}_{\mu^B,\nu^B},    
\end{equation}

where the second sum is taken over tuples of partitions $(\mu^B)_{B\in\pi}$ and $(\nu^B)_{B\in\pi}$ such that $\cup_{B\in\pi}\mu^B=\mu$, $\cup_{B\in\pi}\nu^B=\nu$, and for each $B\in\pi$ we have $|\mu^B|=|\nu^B|=|\lambda_B|$.

\begin{proof}
From equations \eqref{eqtau} and $\eqref{eq d}$ we have 
$$\sum_{n\geq1}\sum_{\lambda,\mu,\nu\vdash n}t^n \frac{d^\lambda_{\mu,\nu}}{z_\lambda(1+b)^{\ell(\lambda)}}p_\lambda q_\mu r_\nu=
\log\left(1+\sum_{n\geq1}\sum_{\lambda,\mu,\nu\vdash n}t^n \frac{c^\lambda_{\mu,\nu}}{z_\lambda(1+b)^{\ell(\lambda)}}p_\lambda q_\mu r_\nu\right).$$
We fix three partitions $\lambda,\mu,\nu$ of the same size.
By developing the logarithm in the previous equation we obtain
$$d^\lambda_{\mu,\nu}=\sum_{s\geq1}\frac{(-1)^{s-1}}{s}\sum_{(\lambda^i,\mu^i,\nu^i)_{1\leq i\leq s}}\frac{z_\lambda}{z_{\lambda^1}...z_{\lambda^s}}\prod_{1\leq i\leq s}c^{\lambda^i}_{\mu^i,\nu^i},$$
where the second sum is taken over $s$-tuples $(\lambda^i,\mu^i,\nu^i)_{1\leq i\leq s}$ such that $\cup_{1\leq i\leq s}\lambda^i=\lambda$, $\cup_{1\leq i\leq s}\mu^i=\mu$ and $\cup_{1\leq i\leq s}\nu^i=\nu$.
We will prove that for every $s\geq 1$,
\begin{equation}\label{Eq s}
  \frac{1}{s!}\sum_{(\lambda^i,\mu^i,\nu^i)_{1\leq i\leq s}}\frac{z_\lambda}{z_{\lambda^1}...z_{\lambda^s}}\prod_{1\leq i\leq s}c^{\lambda^i}_{\mu^i,\nu^i} =\sum_{\pi\in\mathcal{P}_s\left(\llbracket\ell(\lambda)\rrbracket\right)}\sum_{(\mu^B)_{B\in\pi},(\nu^B)_{B\in\pi}}\prod_{B\in\pi}c ^{\lambda_B}_{\mu^B,\nu^B}.
\end{equation}
From a $s$-tuple of partitions $(\lambda^i)_{1\leq i\leq s}$ such that $\cup_i \lambda^i=\lambda$, we can obtain an ordered set-partition  of $\llbracket \ell(\lambda) \rrbracket$ into $s$ sets $B_1$,...,$B_s$ by reordering the parts of $\lambda^1,...,\lambda^s$ to reconstruct $\lambda$.  In this reordering we have $\binom{m_j(\lambda)}{m_j(\lambda^1),...,m_j(\lambda^s)}$ ways to reorder the parts of size $j$. Hence we have 
$$\prod_{1\leq j}\binom{m_j(\lambda)}{m_j(\lambda^1),...,m_j(\lambda^s)}=\frac{z_\lambda}{z_{\lambda^1}...z_{\lambda^s}}$$
ways to obtain an ordered partition of $\llbracket \ell(\lambda)\rrbracket$. Finally, we divide by $s!$, since the first sum in the right hand-side of \Cref{Eq s} is taken over unordered set-partitions. 
\end{proof}
\end{lem}

We deduce the main theorem of this section.
\begin{thm}\label{cor integrality d}
For all partitions $\lambda,\mu,\nu\vdash n\geq 1$, the coefficient $d^\lambda_{\mu,\nu}$ is a polynomial in $b$ with integer coefficients.
\begin{proof}
This is a direct consequence of \cref{lem c-d} and \cref{thm integrality}.
\end{proof}
\end{thm}

\subsection{A generalization for coefficients with \texorpdfstring{$k+2$}{} parameters}\label{ssec k parameters}

In this subsection, we fix $k\geq1$. We consider a generalization of the coefficients $c^\lambda_{\mu,\nu}$ and $h^\lambda_{\mu,\nu}$ as follows; we consider $k+2$ sequences of power-sum variables $\mathbf{p}:=(p_1,p_2...)$ and $\mathbf{q}^{(i)}:=(q^{(i)}_1,q^{(i)}_2...)$ for $0\leq i\leq k$.
For every  $\lambda,\mu^0,...,\mu^k\vdash n\geq1$ we define $c^\lambda_{\mu^0,...,\mu^k}$ and $h^\lambda_{\mu^0,...,\mu^k}$ by 
$$\sum_{n\geq0}t^n\sum_{\theta\vdash n}\frac{1}{j^{(\alpha)}_\theta}J^{(\alpha)}_\theta(\mathbf{p})J^{(\alpha)}_\theta(\mathbf{q}^{(0)})...J^{(\alpha)}_\theta(\mathbf{q}^{(k)})=1+\sum_{n\geq1}t^n\sum_{\lambda,\mu_0,...,\mu_k\vdash n}\frac{c^\lambda_{\mu_0,...,\mu_k}(b)}{z_\lambda(1+b)^{\ell(\lambda)}}p_\lambda q^{(0)}_{\mu_0}...q^{(k)}_{\mu_k},$$
and 
\begin{align*}
&(1+b)\frac{t\partial}{\partial t}\log\left(\sum_{n\geq0}t^n\sum_{\theta\vdash n}\frac{1}{j^{(\alpha)}_\theta}  J^{(\alpha)}_\theta(\mathbf{p})J^{(\alpha)}_\theta(\mathbf{q}^{(0)})...J^{(\alpha)}_\theta(\mathbf{q}^{(k)})\right)\\
&\hspace{8cm}=\sum_{n\geq1}t^n\sum_{\lambda,\mu_0,...,\mu_k\vdash n}\frac{h^\lambda_{\mu_0,...,\mu_k}(b)}{z_\lambda(1+b)^{\ell(\lambda)}}p_\lambda q^{(0)}_{\mu_0}...q^{(k)}_{\mu_k}.    
\end{align*}
It is easy to see that when $k=1$ we recover the definition of the coefficients $h^\lambda_{\mu,\nu}$ given in \Cref{equation h}. These coefficients are related to the enumeration of constellations and are the object of a generalized version of the Matching-Jack conjecture and the $b$-conjecture, see \cite{CD20,BD21}.
The coefficients $c^\lambda_{\mu_0,...,\mu_k}$ satisfy the following multiplicativity property (see \cite[Proposition 6.1]{BD21}).
\begin{prop}
Let $k\geq2$ and $\lambda,\mu^0,...,\mu^k\vdash n\geq1$.
We have
$$c^\lambda_{\mu^0,...,\mu^k}(b)=\sum_{\xi\vdash n}c^\lambda_{\mu^0,...\mu^{k-2},\xi}(b)c^\xi_{\mu^{k-1},\mu^k}(b).$$
\end{prop}

The following corollary is a consequence of the previous proposition and \cref{thm integrality}.
\begin{cor}
For every $k\geq1$ and  $\lambda,\mu^0,...,\mu^k\vdash n\geq 1$, the coefficient $c^\lambda_{\mu^0,...,\mu^k}$ is a polynomial in $b$ with integer coefficients. 
\end{cor}

One can see that the arguments\footnote{We use here a generalized version of \cref{thm polynomiality h} that gives the polynomiality of $h^\lambda_{\mu^0,...,\mu^k}$ for $k>1$, see \cite[Theorem 6.6]{BD21}.} used in \cref{ssec bconj} can be generalized to obtain integrality information about the integrality of coefficients $h^\lambda_{\mu^0,...,\mu^k}$:
\begin{prop}\label{prop gen h}
For every $k\geq1$ and $\lambda,\mu^0,...,\mu^k\vdash n\geq 1$, we have that $\frac{z_\lambda}{n}h^\lambda_{\mu^0,...,\mu^k}$ is an integer polynomial in $b$.
\end{prop}
In particular, since $z_{[n]}=n$, we have the following corollary.
\begin{cor}\label{cor h}
For every $k\geq 1$ and $\mu^0,...,\mu^k\vdash n\geq 1$, the coefficient $h^{[n]}_{\mu^0,...,\mu^k}$ is an integer polynomial in $b$.
\end{cor}
The case of one single part partition in the $b$- and the Matching-Jack conjecture has been considered in previous works and some partial results have been proved in this direction, see \cite{KV16,Dol17,KPV18}. \cref{cor h} establishes the integrality  for the one single part partition case in the "generalized" $b$-conjecture.
This result is known for $k=1$. Indeed, \cite[Theorem 5.10]{CD20} implies that $h^{[n]}_{\mu,\nu,[n],...,[n]}$ is a non-negative integer polynomial in $b$. \cref{prop gen h} is however new for general $k$ and seems to be a step towards the integrality in the $b$-conjecture.

\noindent {\bf Acknowledgements.} The author would like to thank his advisors Valentin Féray and Guillaume Chapuy for many interesting discussions and for reviewing this paper several times.
He also thanks Guillaume Chapuy and Maciej Do\l{}{\k{e}}ga for suggesting to consider the multiplicativity property in order to deduce properties of the coefficients $c^\lambda_{\mu,\nu}$.

\bibliographystyle{halpha}
\bibliography{biblio.bib}
\end{document}